\documentclass{amsart}

\usepackage{latexsym, graphics, psfrag, amscd, amssymb} 

\newtheorem{thm}{Theorem}
\newtheorem{tm}[thm]{Theorem}

\newtheorem{pro}[thm]{Proposition}
\newtheorem{cor}[thm]{Corollary}
\newtheorem{lm}[thm]{Lemma}

\newtheorem{df}[thm]{Definition}

\newcommand{\De}{\Delta}
\newcommand{\Si}{\Sigma}
\newcommand{\Om}{\Omega}
\newcommand{\Ga}{\Gamma}

\newcommand{\al}{\alpha}
\newcommand{\si}{\sigma}
\newcommand{\ep}{\epsilon}
\newcommand{\na}{\nabla}
\newcommand{\la}{\lambda}

\newcommand{\tg}{ {\tilde{g}} }

\newcommand{\tnu}{\tilde{\nu}}

\newcommand{\hg }{ {\hat{g}} }
\newcommand{\hH }{\hat{H}}
\newcommand{\hGa}{\hat{\Gamma}}

\newcommand{\hnu}{\hat{\nu}}

\newcommand{\bc }{\bar{c}}
\newcommand{\be }{\bar{e}}
\newcommand{\bg }{ {\bar{g}} }
\newcommand{\bH }{\bar{H}}
\newcommand{\bnu}{\bar{\nu}}
\newcommand{\bna}{\bar{\na}}

\newcommand{\bOm}{\bar{\Om}}
\newcommand{\bGa}{\bar{\Ga}}
\newcommand{\bth}{\bar{\theta}}
\newcommand{\bsi}{\bar{\si}}

\newcommand{\bR}{\mathbb{R}}

\newcommand{\vD}{\mathbf{D}}
\newcommand{\vX}{\mathbf{X}}
\newcommand{\vY}{\mathbf{Y}}
\newcommand{\vN}{\mathbf{N}}
\newcommand{\vH}{\mathbf{H}}

\newcommand{\loc}{ {\mathrm{loc}} }
\newcommand{\Area}{\mathrm{Area}}
\newcommand{\di}{\mathrm{div}}
\newcommand{\ii}{{\mathrm{in}}}
\newcommand{\ee}{{\mathrm{ex}}}

\newcommand{\inn}[2]{\langle{#1},{#2}\rangle} 
\newcommand{\norm}[1]{\parallel{#1}\parallel}

\newcommand{\pa}{\partial}
\newcommand{\D}{\bar{\vD}}
\newcommand{\vvD}{\check{\vD}}
\newcommand{\dep}{\delta,p}

\begin{document}

\title{Positivity of quasi-local mass II}

\author{Chiu-Chu Melissa Liu}
\address{Department of Mathematics\\ Harvard University\\ Cambridge, MA
02138, USA}
\email{ccliu@math.harvard.edu}

\author{Shing-Tung Yau}
\address{Department of Mathematics\\ Harvard University\\ Cambridge, MA
02138, USA}
\email{yau@math.harvard.edu}

\begin{abstract}
We prove the following stronger version of the positivity of quasi-local mass
stated in \cite{Li-Ya}: the quasi-local energy (mass) of each connected component
of the boundary of a compact spacelike hypersurface which satisfies the local
energy condition is strictly positive unless the spacetime is flat along the
spacelike hypersurface and the boundary of the spacelike hypersurface is
connected. 
\end{abstract}

\maketitle

\section{Introduction} \label{sec:introduction}

A spacetime is a four manifold with a pseudo-metric of signature $(+,+,+,-)$.
A hypersurface or a 2-surface in a spacetime is spacelike if the
induced metric is positive definite. Quasi-local energy-momentum
vector is a vector in $\bR^{3,1}$ associated to a spacelike 2-surface
which depends on its first and second fundamental forms, and 
the connection on its normal bundle in the spacetime. The time
component of the four vector is called quasi-local mass. Similar to
\cite{Ea, Ch-Ya}, we require the quasi-local energy-momentum vector 
to satisfy the following properties.
\begin{enumerate}
\item It should be zero for the flat spacetime.
\item The quasi-local mass should be equivalent to the standard definition if the spacetime is
      spherically symmetric and quasi-local mass is evaluated on the spheres 
      \cite{Ca-Mc}. (We say two masses $m_1$ and $m_2$ are
      equivalent if there is a universal constant $c>0$ such
      that $c^{-1}m_1\leq m_2\leq c m_1$.)
      In particular, for the centered spheres in the Schwarzschild spacetime,  
      the quasi-local mass should be equivalent to the standard mass.
\item For an asymptotically flat slice, the quasi-local mass of the coordinate
      sphere should be asymptotic to the ADM energy-momentum vector.
\item For an asymptotically null slice, the quasi-local mass of the coordinate
      sphere should be asymptotic to the Bondi energy-momentum vector.
\item For an apparent horizon $\Si$, the quasi-local mass should be no
      less than a (universal) constant multiple of the irreducible
      mass which is $\sqrt{\Area(\Si)/16\pi}$.
\item The quasi-local energy-momentum vector should be non-spacelike
      and the quasi-local mass should be non-negative.
\end{enumerate}

Our definition of quasi-local mass \cite{Li-Ya} arises naturally from
calculations in the second author's work \cite{Ya} 
on blackholes, and is strongly motivated by our ability to prove its
positivity. After the second author proposed our definition, we were
informed of the existence of much earlier works by Brown-York
\cite{Br-Yo1, Br-Yo2} and others \cite{La, Ki, Ep}.
The main goal of this paper is to provide 
a complete proof of a stronger version of the positivity stated in \cite{Li-Ya}. 

The rest of the paper is organized as follows.
In Section \ref{sec:statement}, we recall our definition of
quasi-local mass, discuss its properties, and
state the main result (positivity of quasi-local mass). 
In Section \ref{sec:shi-tam}, we describe
Shi and Tam's proof of positivity in the Riemannian case. 
In Section \ref{sec:proof}, we prove the main result.

\bigskip

\paragraph{\bf Acknowledgments.} 
We wish to thank Shi and Tam for their questions on \cite{Li-Ya} which
motivated us to write down this proof in detail. We wish to thank 
Mu-Tao Wang for proofreading the first draft. The first author
also wishes to thank Xiao Zhang for helpful discussion.
The second author is supported in part by the National Science 
Foundation under Grant No. DMS-0306600. 

\section{Definition of quasi-local mass and its properties}\label{sec:statement}
Let $\Si$ be a spacelike 2-surface in a spacetime $N$.
At each point of $\Si$, choose two null normals $l, n$ such that
$\inn{l}{n}=-1$. 
Any other choice $(l',n')$ is related to $(l,n)$ by
$l'=\la l, n'=\la^{-1}n$ or $l'=\la n, n'=\la^{-1}l$
for some function $\la:\Si\to \bR\setminus\{0\}$. We
denote the mean curvature with respect to $l$ and $n$ by 
\begin{equation}\label{eqn:null}
2\rho=-\inn{\na_1 e_1+\na_2 e_2}{l},\  \
-2\mu=-\inn{\na_1 e_1+\na_2 e_2}{n}
\end{equation}
respectively, where $\{e_1,e_2\}$ is a local orthonormal frame of $\Si$.
The definitions of $\rho$ and $\mu$ depend on choice of $(l,n)$, but their product
$\rho\mu$ is independent of choice of $(l,n)$.  More intrinsically,
$$
8\rho\mu=\inn{\vH}{\vH}
$$
where $\vH$ is the mean curvature vector of $\Si$ in $N$.
We assume that $\rho\mu>0$, or equivalently,
the mean curvature vector $\vH$ of $\Si$
in $N$ is spacelike.

Suppose that $\Si$ has positive Gaussian curvature so that 
$\Si$ is topologically a 2-sphere.  By Weyl's embedding
theorem, $\Si$ can be isometrically embedded into the Euclidean
space $\bR^3$ so that the second fundamental form $(H_0)_{ab}$
is positive definite. The embedding $\Si\subset\bR^3$ is unique
up to an isometry of $\bR^3$, so $(H_0)_{ab}$ is determined by
the metric on $\Si$. Let $\rho_0, \mu_0$ be the mean curvatures
with respect to null normals $l_0,n_0$  of the embedding
$\Si\subset\bR^3\subset\bR^{3,1}$,
with the normalization $\inn{l_0}{n_0}=-1$. Then
$$
8\rho_0\mu_0=H_0^2
$$
where $H_0>0$ is the trace of $(H_0)_{ab}$.
Define the {\em quasi-local mass} of $\Si$ to be
\begin{equation}\label{eqn:energy}
E(\Si)=\frac{1}{8\pi G}\int_\Si (\sqrt{8\rho_0\mu_0}-\sqrt{8\rho\mu})
=\frac{1}{8\pi G}\int_\Si(H_0-\sqrt{8\rho\mu})
\end{equation}
See \cite{Sz} for other definitions of quasi-local mass (energy).

Recall (1)--(6) in Section \ref{sec:introduction}.
In \cite{MST}, Murchadha, Szabados, and Tod gave examples of
$\Si\subset \bR^{3,1}$ but $E(\Si)>0$, so $E(\Si)$ does not satisfy (1). For (2),
recall that the Schwarzschild spacetime metric on $\bR^4$ is given by
$$
g=-\left(1-\frac{2M}{r}\right)dt^2+\left(1-\frac{2M}{r}\right)^{-1} dr^2
+r^2(d\theta^2+\sin^2\theta d\phi^2),\ \  r>2M
$$
where $r,\theta,\phi$ are the spherical coordinates on $\bR^3$. Let
$S_a\subset (\bR^4,g)$ be the round sphere defined by $t=0, r=a$,
and let $m(r)=E(S_r)$. Then
$$
m(r)=r(1-\sqrt{1-\frac{2M}{r} }).
$$
Note $m(r)$ is decreasing (for $r\geq 2M$), 
$m(2M)=2M$, and $m(\infty)=M$, which is consistent with (2).
For (3), (4), Epp discussed the
spatial
and future null infinity limits of a large sphere in asymptotically
flat spacetimes, but cannot conclude that $E(\Si)$ satisfies (3), (4) in general.
For (5), on an apparent horizon $\Si$ we have $\rho\mu=0$, so
$$
E(\Si)=\frac{1}{8\pi}\int_{\Si} H_0 \geq
\sqrt{\frac{\mathrm{Area(\Si)}}{4\pi} }
$$
by Minkowski inequality of convex bodies \cite{Mi}. Therefore, $E(\Si)$ satisfies (5).
By the main result of this paper, $E(\Si)$ is nonnegative as required
in (6) and is strictly positive when the spacetime is not flat along $\Si$. 

We now give precise statement of the main result.
Let $\Om$ be a compact spacelike hypersurface in a time orientable
four dimensional spacetime $N$. Let $g_{ij}$ denote the 
induced metric on $\Om$, and let $p_{ij}$ denote the second
fundamental form of $\Om$ in $N$. The local mass density $\mu$ and 
the local current density $J^i$ on $\Om$ are related to $g_{ij}$ and 
$p_{ij}$ by the constraint equations
\begin{eqnarray} \label{eqn:constriant}
\mu&=&\frac{1}{2}\left(R-\sum_{i,j}p^{ij}p_{ij} + (\sum_i p^i_i)^2\right)\\
J^i&=&\sum_j D_j\left(p^{ij}-(\sum_k p_k^k)g^{ij}\right),
\end{eqnarray}
where $R$ is the scalar curvature of the metric $g_{ij}$.
In this paper, we prove the following stronger version of
the positivity stated in \cite{Li-Ya}.
\begin{tm}[positivity of quasi-local mass]\label{thm:positive}
Let $\Om,\mu,J$ be as above.
We assume that $\mu$ and $J^i$ satisfies the local energy condition
\begin{equation}\label{eqn:local}
\mu\geq\sqrt{J^i J_i}
\end{equation}
and the boundary $\pa\Om$ has finitely many connected components 
$\Si^1,\ldots,\Si^\ell$, each of which has positive Gaussian curvature 
and has spacelike mean curvature vector in $N$.
Let $E(\Si^\al)$ be defined as in (\ref{eqn:energy}). 
Then $E(\Si^\al)\geq 0$ for $\al=1,\ldots,\ell$. 
Moreover, if $E(\Si^\al)=0$ for some $\al$, then $M$ is flat spacetime 
along $\Om$, $\pa\Om$ is connected and will be embedded into 
$\bR^3\subset\bR^{3,1}$ by the well-known Weyl embedding theorem. 
\end{tm}

\section{The Riemannian case}\label{sec:shi-tam}
When the second fundamental form of $\Om$ in $N$ vanishes, 
the local energy condition (\ref{eqn:local}) reduces to $R\geq 0$ and
the condition $\rho\mu>0$ reduces to $H>0$, where $H$ is the mean
curvature of the spacelike 2-surface in $\Om$ with respect to the outward
unit normal. Shi and Tam proved positivity of quasi-local mass in this case. 

\begin{tm}[{\cite[Theorem 1]{Sh-Ta}}]\label{thm:shi-tam}
Let $(\Om^3,g)$ be a compact manifold of dimension three with smooth
boundary and with nonnegative scalar curvature. Suppose $\pa\Om$ has 
finitely many connected components $\Si^\al$ so that each connected 
component has positive Gaussian curvature and positive mean curvature
$H$ with respect to the unit outward normal. Then for each boundary 
component $\Si^\al$,
\begin{equation}\label{eqn:mean}
\int_{\Si^\al}H d\si\leq\int_{\Si^\al}H_0^\al d\si
\end{equation}
where $H_0^\al$ is the mean curvature of $\Si^\al$ with respect to the outward
normal when it is isometrically embedded in $\bR^3$, $d\si$ is the volume 
form on $\Si^\al$ induced from $g$. Moreover, if equality holds in 
(\ref{eqn:mean}) for some $\Si^\al$, then $\pa\Om$ has only one connected
component and $\Om$ is a domain in $\bR^3$.
\end{tm} 

We now briefly describe Shi and Tam's proof of Theorem \ref{thm:shi-tam}.
From now on, all the mean curvatures are defined with respect to
the outward unit normal. 

Let $\Si^\al$ be a connected component of $\pa\Om$. By hypothesis
of Theorem \ref{thm:shi-tam},
it has positive Gaussian curvature, so it can be isometrically
embedded to $\bR^3$ by the well-known Weyl embedding theorem.
Moreover, the embedding is unique up to an isometry of $\bR^3$.
Let $\Si^\al_0\subset\bR^3$ be the image of such an embedding. 
Then $\Si^\al_0$ is a strictly convex hypersurface diffeomorphic to $S^2$. 

Let $\vX$ be the position vector of a point on $\Si^\al_0$, and let
$\vN$ be the unit outward normal of $\Si^\al_0$ at $\vX$.
Let $\Si^\al_r$ be the surface described by $\vY=\vX+r\vN$,
with $r\geq 0$. Let $D^\al$ be the region of $\bR^3$ outside $\Si^\al_0$, and
let $E^\al=D^\al\cup \Si^\al_0$ be the closure of $D^\al$ in
$\bR^3$. Then $E^\al$ can be represented by
\[
(\Si^\al\times [0,\infty), g^0=dr^2 + g_r),
\]
where $g_r$ is the induced metric on $\Si^\al_r$, $g^0=dr^2+g_r$ is the standard
Euclidean metric on $E^\al\subset\bR^3$. Note that 
$g_r=(a+r)^2(d\theta^2+\sin^2\theta d\phi^2)$ if $\Si^\al_0$ is a
round sphere of radius $a>0$.

Consider a Riemannian metric on $E^\al$ of the form
\begin{equation}
g=h^2 dr^2 + g_r.
\end{equation}
where $h$ is a smooth positive function. This is a special case of
Bartnik's construction in \cite{Ba2}. 
Note that $g$ and $g^0$ induce the same metric on each $\Si^\al_r$. The mean
curvature $H$ and $H_0$ of $\Si^\al_r$ with respect to 
$g$ and $g^0$ are related by
\[
H=h^{-1}H^0. 
\]
Note that $H_0(x,0)=H_0^\al(x)$ for $x\in \Si^\al_0\cong\Si^\al$.
The scalar curvature $R$ of $g$ is given by
\[
2H_0\frac{\pa h}{\pa r}=2h^2\Delta_r h+ (h-h^3)R^r+ h^3 R,
\]
where $R^r$ is the scalar curvature of $\Si^\al_r$, and $\Delta_r$
is the Laplacian operator on $\Si^\al_r$.
So a solution to the parabolic partial 
differential equation
\begin{equation}\label{eqn:shi-tam-foliation}
2H_0\frac{\pa h}{\pa r} =2h^2\Delta_r h+(h-h^3)R^r 
\end{equation}
on $E^\al\cong \Si^\al\times[0,\infty)$
with the initial condition
\begin{equation}\label{eqn:shi-tam-initial}
h(x,0)=\frac{H_0^\al}{H}.
\end{equation}
defines a metric on $E^\al$ such that the scalar curvature $R=0$ and 
the mean curvature of $\Si^\al_0$ coincides with the restriction of $H$ 
to $\Si^\al\cong \Si^\al_0$.

Let $\rho:\bR^3\to [0,\infty)$ be the distance function 
to the origin (in Euclidean metric). We may assume that
the origin is enclosed by $\Si^\al_0$ so that $\rho \geq a$
for some constant $a>0$. Shi and Tam showed that 
\cite[Theorem 2.1]{Sh-Ta}:
\begin{tm}\label{thm:unique}
The equation (\ref{eqn:shi-tam-foliation}) with the initial condition 
(\ref{eqn:shi-tam-initial}) has a unique solution such that
\begin{enumerate}
\item[(a)] $h=1+ m_o \rho^{-1} + \kappa$, where $m_0$ is a constant
and the function $\kappa$ satisfies 
$$
|\kappa|=O(\rho^{-2}), \ \ \ 
|\nabla_0 \kappa|=O(\rho^{-3}),
$$
where $\nabla_0$ is Levi-Civita connection of the Euclidean
metric on $\bR^3$.
\item[(b)] The metric $g^\al=h^2dr^2 + g_r$ on $E^\al$ 
is asymptotically flat in the sense that 
\begin{equation}\label{eqn:asymflat}
|g^\al_{ij}-\delta_{ij}| + \rho |\nabla_0 g^\al_{ij}| 
+\rho^2 |\nabla^2_0 g^\al_{ij}|\leq C\rho^{-1}.
\end{equation}
with zero scalar curvature.
\item[(c)] The ADM mass of $(E^\al,g^\al)$ is given by
$$
m^\al_\infty=\lim_{r\to \infty}m^\al(r),
$$ 
where
$$
m^\al(r)=\frac{1}{8\pi G}\int_{\Si^\al_r} (H_0-H)d\si_r,
$$
and $d\si_r$ is the volume form of $\Si^\al_r$.
\end{enumerate}
\end{tm}

Let $m^\al(r)$ be defined as in (c) of Theorem \ref{thm:unique}.
It is computed in the proof of \cite[Lemma 4.2]{Sh-Ta} that
\begin{equation} \label{eqn:monotone}
\frac{d m^\al}{dr}(r)=\frac{-1}{16\pi G}\int_{\Si^\al_r} R^r u^{-1}(1-u)^2 \leq 0.
\end{equation}

Shi and Tam glued $(E^\al,g^\al)$ to $(\Om,g)$ along
$\Si^\al$ to obtain a complete noncompact three manifold $M$
with a continuous Riemannian metric $\tg$ such that
\begin{enumerate}
\item $\tg$ is smooth on $M\setminus\Om$ and $\bOm$, and is Lipschitz
      near $\pa\Om$.
\item The mean curvatures of $\Si^\al$ with respect to
      $g=\tg|_{\Om}$ and $g^\al=\tg|_{E^\al}$ are the same for each
      $\al$.
\item Each end $E^\al$ of $M$ is asymptotically Euclidean in the 
      sense of (\ref{eqn:asymflat}).
\item The scalar curvature $R$ of $M\setminus\pa\Om$ is nonnegative and
      is in $L^1(M)$.
\end{enumerate}
Using Witten's argument \cite{Wi,Pa-Ta}, Shi and Tam proved that the positive mass theorem holds for 
such a metric, so the ADM mass 
\[
m_\infty^\al=\lim_{r\to\infty}m^\al(r)
\]
is nonnegative for each end $E^\al$, and $m_\infty^\al$ vanishes for some $\al$
if and only if $M$ has only one end and $M$ is flat. This together
with the monotonicity (\ref{eqn:monotone}) of $m^\al(r)$ gives 
Theorem \ref{thm:shi-tam}, since
\[
m^\al(0)=\frac{1}{8\pi G}\int_{\Si^\al} (H_0^\al - H)d\si.
\]

\section{Proof of Theorem \ref{thm:positive}}\label{sec:proof}
\subsection{Outline of proof}\label{sec:general}

Let $(\Omega,g_{ij},p_{ij})$ and $\Si^1,\ldots,\Si^\ell$ be as
in Section \ref{sec:statement}. We first deform the metric $g_{ij}$
on $\Omega$ by a procedure used by Schoen and the second
author in \cite{Sc-Ya} and also by the second author in \cite{Ya}.
This procedure consists of two steps. The first step is to deform
$g_{ij}$ to a new metric   
$$
\bg_{ij}= g_{ij} + f_i f_j
$$ 
where $f$ is a solution to Jang's equation on $\Omega$ such that
$f|_{\pa\Om}=0$. The metric $\bg_{ij}$ coincides with $g_{ij}$ when restricted to
$\pa\Om$, and its scalar curvature
$\bar{R}$ satisfies
\begin{equation}\label{eqn:RX}
\bar{R}\geq 2|X|^2-2\di X
\end{equation}
for some vector field $X$ on $\Omega$. The equality holds only
if $p_{ij}=h_{ij}$, where $h_{ij}$ is the second fundamental form 
of the isometric embedding of  $(\Om,\bg_{ij})$ into
$$
(\Om\times \bR, g_{ij} dx^i dx^j +dt^2)
$$ 
as the graph of $f$. The second step is to deform 
$\bg_{ij}$ conformally to a metric with zero scalar
curvature. The inequality (\ref{eqn:RX}) implies
that there is a unique metric $\hg_{ij}$ in the conformal 
class of $\bg_{ij}$ which has zero scalar curvature and coincides with
$\bg_{ij}$ on $\pa\Om$. 

After the above reduction, we cannot apply Theorem \ref{thm:shi-tam}
directly to $\hg_{ij}$ because the mean curvature $\hH$ of $\pa\Om$
with respect to $\hg_{ij}$ is not necessarily positive (this point was
overlooked in \cite{Li-Ya}). Instead, we have
\begin{equation}\label{eqn:hHbH}
\int_{\pa\Om}\hH \geq \int_{\pa\Om} (\bH-\inn{X}{\bnu})
\end{equation}
where $\bH$ and $\bnu$ are the mean curvature  and outward unit normal
of $\pa\Om$ with respect to $\bg_{ij}$, and the equality holds iff 
$\hg=\bg$ and $X=0$ (so $p_{ij}=h_{ij}$). It was shown in \cite{Ya} that
\begin{equation}\label{eqn:HXHP}
\bH-\inn{X}{\bnu}\geq \sqrt{H^2-P^2}=\sqrt{8\rho\mu}
\end{equation}
where $P$ is the trace of the restriction of $p_{ij}$ to $\pa\Om$.
In particular, $\bH-\inn{X}{\bnu}$ is positive.

Let $E^\alpha$ and $H_0^\alpha$ be defined as
in Section \ref{sec:shi-tam}. Shi and Tam's proof of
Theorem \ref{thm:unique} shows that one can solve (\ref{eqn:shi-tam-foliation}) 
on $E^\alpha$ with the initial condition
\begin{equation}
h(x,0)=\frac{H_0^\alpha}{\bH-\inn{X}{\bnu}}
\end{equation}
and obtain a scalar flat, asymptotically flat metric $g^\alpha$ on
the end $E^\alpha$. Gluing $(E^\alpha, g^\alpha)$ to
$(\Omega, \hg)$ along $\Si^\alpha$, we obtain a complete noncompact
three manifold $M$ with a Lipschitz continuous Riemannian metric
$\tg$. On $M\setminus \pa\Om$, $\tg$ is smooth and has zero scalar
curvature. However, the mean curvature of $\Si^\alpha$ with respect
to $\hg=\tg|_\Om$ and $g^\alpha=g|_{E^\al}$ are not necessarily the
same. This causes the following problem which is absent in the case 
considered by Shi and Tam: the zeroth order term of the Dirac operator 
can be discontinuous along $\pa\Om$, so there is an extra term when
we integrate the Weitzenb\"{o}ck-Lichnerowicz formula. To prove the positive
mass theorem for $(M,\tg)$ (Theorem \ref{thm:tgmass}), we derive an inequality (Proposition
\ref{thm:key}) as a substitute of the integral form of
Weitzenb\"{o}ck-Lichnerowicz formula for smooth metrics. 
 
Let $m^\alpha_\infty$ and $m^\alpha(r)$ be defined by 
$g^\alpha$ as in Section \ref{sec:shi-tam}. The monotonicity
(\ref{eqn:monotone}) of $m^\alpha(r)$  and (\ref{eqn:HXHP}) imply
\begin{equation}
m^\al_\infty\leq \frac{1}{8\pi G}\int_{\Si^\alpha}(H_0^\alpha -
(\bH-\inn{X}{\bnu})) \leq \frac{1}{8\pi
  G}\int_{\Si^\alpha}(H_0^\alpha-\sqrt{8\rho\mu})
=E(\Si^\alpha).
\end{equation}
The positive mass theorem for $(M,\tg)$ says that $m^\al_\infty\geq 0$ for 
$\al=1,\ldots,\ell$, and $m^\al_\infty=0$ for some $\al$ iff $\ell=1$ and
$(M,\tg)$ is the Euclidean space $\bR^3$. So $E(\Si^\al)\geq 0$ for 
$\al=1,\ldots,\ell$. If $E(\Si^\al)=0$ for some $\al$, we must have
$m^\al_\infty=0$ and $\hg=\bg$, so $(\Omega,\bg)=(\Omega,\tg)$ is a domain
$\Omega_0\subset \bR^3$. In this case, $(\Omega,g)$ (at least the
part away from apparent horizons) can be
isometrically embedded in $\bR^{3,1}=(\bR^3\times \bR, \sum
dx_i^2-dt^2)$ as 
a graph 
$$
\{ (x,f(x))\mid x\in \Omega_0)\}
$$
with second fundamental form $p_{ij}$, where $f$ is a smooth function 
on $\Om_0$ which vanishes on $\pa\Om_0$.

\subsection{Jang's equation with Dirichlet boundary condition}
\label{sec:jang}

As in \cite{Sc-Ya}, we consider following equation proposed by 
Jang \cite{Ja} on $\Omega$:
\begin{equation}\label{eqn:jang}
\sum_{i,j=1}^3\left(g^{ij}-\frac{f^i f^j}{1+|\na f|^2}\right)
\left(\frac{f_{ij}}{\sqrt{1+|\na f|^2}}-p_{ij}\right)=0.
\end{equation}
As in \cite{Ya}, we consider solutions to (\ref{eqn:jang}) 
with the Dirichlet boundary condition
\begin{equation}\label{eqn:zero}
f|_{\pa\Omega}\equiv 0.
\end{equation}

Most of the estimates were made in \cite{Sc-Ya}. To solve the boundary
value problem, the second author constructed a barrier in
\cite{Ya} and concluded that there exists a solution to (\ref{eqn:jang}) with
boundary value (\ref{eqn:zero}) when
$(\Omega,g_{ij},p_{ij})$ has no {\em apparent horizon}.

\begin{df}
Let $(\Omega, g_{ij},p_{ij})$ be an initial data set.
Given a smooth compact surface $S$ embedded
in $\Omega$, let $H_s$ be the mean curvature of $S$ with respect to the
outward unit normal vector, and let $P_s$ be the trace of the restriction
of $p_{ij}$ to $S$. A smooth 2-sphere $S$ embedded in $\Om$ is an {\em apparent horizon}
of the initial data $(\Omega,g_{ij},p_{ij})$ if $H_s+P_s=0$ or $H_s-P_s=0$. 
\end{df}

We first assume that $(\Omega,g_{ij},p_{ij})$ has no apparent horizon
so that there exists a solution $f$ to the Jang's equation
(\ref{eqn:jang}) on $\Omega$ such that $f|_\Om=0$.
The induced metric of the graph $\Om_f\cong \Om$ of $f$ in 
$(\Om\times\bR, g_{ij}dx^i dx^j + dt^2)$
is 
\[
\bg_{ij}=g_{ij}+f_i f_j
\]
which can be viewed as a deformation of the metric $g_{ij}$ on $\Om$.
Note that the new metric $\bg$ coincides with
the old metric $g$ when restricted to  $\pa\Om$.

We now introduce some notation.
Let $\be_4$ be the downward unit normal to $\Om_f$ in $\Om\times\bR$,
and let $\be_1,\be_2,\be_3$ be a local orthonormal frame of $\Om$. We 
define $h_{i4}$ by
\[
\bna_4 \be_4 = h_{i4}\be_i.
\]
where $\bna$ denotes the Levi-Civita connection of
the metric $g_{ij}dx^i dx^j + dt^2$ on $\Om\times\bR$. 
Let $h_{ij}=\inn{\be_i}{\bar{\na}_j \be_4}$ be the 
second fundamental form of $\Om_f$ in $\Om\times\bR$. 
Let $\bar{R}$ be the scalar curvature of $\bg$, and extend
$p_{ij}$, $\mu$, $J^i$ parallely along the $\bR$ factor. 
The following inequality was derived in \cite{Sc-Ya}:
\begin{equation}\label{eqn:hp}
       2(\mu-|J|)\leq \bar{R}-\sum_{i,j}(h_{ij}-p_{ij})^2
      -2(h_{i4}-p_{i4})^2+ 2\sum_i D_i(h_{i4}-p_{i4}),
\end{equation}
where $D_i$ denotes the covariant derivative of $\bg$.
In particular,
\begin{equation}\label{eqn:RXX}
 \bar{R}\geq 2|X|^2-2\mathrm{div}X,
\end{equation}
where $X=\sum(h_{i4}-p_{i4})e_i$, and the divergence is defined by
$\bg$. By (\ref{eqn:hp}), the inequality (\ref{eqn:RXX}) is 
an equality only if $p_{ij}=h_{ij}$.

In general, the solution $f$ and the metric $\bg$ 
are defined on $\Om'$, the complement of union of apparent horizons,
but one can extend $\bg$ to a metric on $\Om''$ which is obtained by 
adding a point on each end of $\Om'$. See \cite{Sc-Ya} for details. 

\subsection{Scalar flat metric on $\Omega$} \label{sec:interior}
We shall prove the following:
\begin{pro}\label{thm:conformal}
Let $(\Om,\bg)$ be a compact Riemannian manifold of dimension three
with smooth boundary. Suppose that the scalar curvature $\bar{R}$
of $\bg$ satisfies
\[
\bar{R}\geq  c|X|^2-2\mathrm{div}X,
\]
for some constant $c > \frac{1}{2}$ and some smooth vector field
$X$  on $\Om$. Then there is a unique metric $\hg_{ij}$
on $\Om$ such that
\begin{enumerate}
\item The metric $\hg_{ij}$ is conformal to $\bg_{ij}$.
\item The scalar curvature of $\hg_{ij}$ is zero.
\item The metric $\hg_{ij}$ coincides with $\bg_{ij}$ on $\pa\Om$.
\item Let $\bH$ and $\hH$ denote the mean curvatures
with respect the metric $\bg$ and $\hg$, respectively, and
let $\bnu$ denote the outward unit normal of $\pa\Om$ in $(\Om,\bg)$. 
Then
\[
\int_{\pa\Om}\hH \geq \int_{\pa\Om}(\bH-\inn{X}{\bnu}),
\]
where the equality holds if and only of $\bar{R}=0$, $X=0$, and $\hg_{ij}=\bg_{ij}$. 
\end{enumerate}
\end{pro}

\paragraph{Proof.} 
In this proof, the Laplacian, gradient, divergence, and all
the norms are defined by the metric $\bg_{ij}$.

Any metric $\hg_{ij}$ conformal to $\bg_{ij}$ can be written as 
$\hg_{ij}=u^4 \bg_{ij}$,  where $h$ a positive smooth function on $\Om$.
The metric $\hg_{ij}$ satisfies (1) and (2) in Proposition
\ref{thm:conformal} if and only if $v=u-1$ is a solution to 
\begin{equation}\label{shift}
\left\{\begin{array}{cl}
\De v-\frac{1}{8}\bar{R}v=\frac{1}{8}\bar{R}& \textup{ on }\Om \\
v=0 & \textup{ on }\pa\Om
\end{array}\right.
\end{equation}

We first show that (\ref{shift}) has a unique solution.
Let $f$ be a solution to 
\begin{equation}\label{homo}
\left\{\begin{array}{cl}
\De f-\frac{1}{8}\bar{R}f=0& \textup{ on }\Om \\
f=0 & \textup{ on }\pa\Om
\end{array}\right.
\end{equation}
Then
\begin{eqnarray*}
0&=&\int_\Om f\left(-\De f+\frac{\bar{R}}{8}f\right)
=\int_\Om\left(|\na f|^2+\frac{\bar{R}}{8}f^2\right) \\
&\geq &\int_\Om\left(|\na f|^2 +\frac{c}{8}|X|^2 f^2
      -\frac{1}{4}(\mathrm{div}X)f^2\right)\\
&=&\int_\Om \left(|\na f|^2+ \frac{1}{2}X(f)f
   +\frac{c}{8}|X|^2 f^2\right)\\
&\geq&\int_\Om \left(|\na f|^2 -\frac{1}{2}|fX||\na f|
       +\frac{c}{8}|fX|^2 \right)\\
&=& \int_\Om \left(\left(\frac{1}{\sqrt{2c}}|\na f|
   -\frac{\sqrt{2c}}{4}|fX|\right)^2+
       (1-\frac{1}{2c})|\na f|^2\right)\\
&\geq&0.
\end{eqnarray*}

Note that $1-\frac{1}{2c}>0$, so $\na f\equiv 0$, which implies
$f\equiv 0$ since $f$ vanishes on $\pa\Om$. 
Therefore, zero is the only solution to (\ref{homo}), 
and (\ref{shift}) has a unique solution.
Let $v$ be the unique solution to (\ref{shift}). Then
$v$ is smooth. 

We next show that $u=v+1$ is positive. 
Note that $u$ satisfies 
\begin{equation}\label{scalar}
\De u-\frac{1}{8}\bar{R}u=0
\end{equation}
on $\Om$. Assume that $\Om_-=\{ x\in\Om\mid u(x)<0\}$ is nonempty.
Then $\pa\Om_-\cap \pa\Om=\emptyset$, and
\[
\left\{\begin{array}{cl}
\De u-\frac{1}{8}\bar{R}u=0& \textup{ on }\Om_- \\
u=0 & \textup{ on }\pa\Om_-
\end{array}\right.
\]
which implies that $u\equiv 0$ on $\Om_-$, a contradiction.
So $\Om_-$ must be empty, or equivalently, $u$ is nonnegative.
Since $u=1$ on $\pa\Om$, the positivity of $u$ follows from 
the Harnack inequality for nonnegative solutions to (\ref{scalar}).

Finally, we check that the metric $\hg_{ij}=u^4 \bg_{ij}$
satisfies (4) in Proposition \ref{thm:conformal}.
We first compute $\hH$ in terms of $\bH$ and $u$.
We may choose coordinates $(x^1,x^2,x^3)$ near a point
$p\in\pa\Om$ such that $(x^1,x^2)$ are the normal 
coordinates of $\pa\Om$ centered at $p$, and 
$\frac{\pa}{\pa x^3}=\bnu$.
Let $\bGa_{ij}^k$, $\hGa_{ij}^k$ denote the
Christoffel symbols of $\bg$, $\hg$, respectively.
Then for $i,j=1,2$, at $p$ we have 
\begin{eqnarray*}
-\hat{h}_{ij}=\hGa_{ij}^3
&=& \frac{1}{2}\left(\frac{\pa}{\pa x^i}\hg_{i3}
                +\frac{\pa}{\pa x^j}\hg_{j3}
                -\frac{\pa}{\pa x^3}\hg_{ij}\right)\\
&=& \frac{1}{2}\left( \frac{\pa}{\pa x^i}(u^4 \bg_{i3})
                     +\frac{\pa}{\pa x^j}(u^4 \bg_{j3})
                     -\frac{\pa}{\pa x^3}(u^4 \bg_{ij})\right)\\
&=& \frac{1}{2}\left(  \frac{\pa}{\pa x^i}\bg_{i3}
                      +\frac{\pa}{\pa x^j}\bg_{j3}
                      -\frac{\pa}{\pa x^3}\bg_{ij}
                      -4\frac{\pa u}{\pa x^3}\delta_{ij}\right)\\
&=& \bGa_{ij}^3-2\frac{\pa u}{\pa x^3}\delta_{ij}\\
&=&-\bar{h}_{ij}-2\frac{\pa u}{\pa x^3}\delta_{ij}
\end{eqnarray*}
So $\hH=\bH+4\bnu(u)$.

\begin{eqnarray*}
\int_{\pa\Om} \bnu(u) &=& \int_\Om \mathrm{div}(u\na u)
 = \int_\Om(|\na u|^2 + u\De u)
 = \int_\Om\left(|\na u|^2 +\frac{\bar{R}}{8}u^2\right)\\
&\geq& \int_\Om\left(|\na u|^2 + \frac{c}{8}|X|^2 u^2
      - \frac{1}{4}\mathrm{div} X u^2\right)\\
&=&  \int_\Om \left(|\na u|^2 + \frac{c}{8}|X|^2 u^2
      + \frac{1}{2} X(u)u \right)
      -\frac{1}{4}\int_{\pa\Om}\inn{X}{\bnu} \\  
&\geq& \int_\Om \left(\left(\frac{1}{\sqrt{2c} }|\na u|
   -\frac{\sqrt{2c}}{4}|uX|\right)^2+(1-\frac{1}{2c})|\na u|^2\right)
 -\frac{1}{4}\int_{\pa\Om}\inn{X}{\bnu}\\
&\geq& -\frac{1}{4}\int_{\pa\Om}\inn{X}{\bnu} 
\end{eqnarray*}

Therefore, 
\[
\int_{\pa\Om}\hH \geq \int_{\pa\Om}(\bH-\inn{X}{\bnu}),
\]
and the equality holds if and only if 
\[
u\equiv 1,\ X=0,\ \bar{R}=0.\ \ \Box
\]

\subsection{Scalar flat metrics on the ends}\label{sec:end}

\begin{lm}\label{thm:XHP}
\begin{equation}\label{eqn:hxhp}
\bH-\inn{X}{\bnu}\geq \sqrt{H^2-P^2}
\end{equation}
\end{lm}

\paragraph{Proof}
Let $\{ \be_1,\be_2,\be_3,\be_4\}$ be a local orthonormal frame of 
$\Om\times \bR$ along the graph $\Om_f$ so that $\be_1,\be_2$ is tangent
to $\pa\Om$ and $\be_3=\bnu$. Let $w$ be the outward unit normal
of $\pa\Om_0$ in $\Om_0$, the graph of the zero function. 
It was computed in \cite[Section 5]{Ya} that 
\begin{equation} \label{eqn:black}
\bH-\inn{X}{\bnu}=-\frac{\inn{\be_4}{w}}{\inn{\be_3}{w}}P
+\frac{1}{\inn{e_3}{w}}H. 
\end{equation}
Recall that $H>0$, so (\ref{eqn:hxhp}) is equivalent to 
\[
\left( -\inn{\be_4}{w}P+ H \right)^2\geq 
\inn{\be_3}{w}^2 
\left( H^2-P^2 \right)
\]
which is equivalent to
\[ 
(\inn{e_4}{w}^2+\inn{e_3}{w}^2)P^2
-2\inn{e_4}{w}PH + (1-\inn{e_3}{w}^2) H^2
\geq 0
\]
But $|w|^2=1=\inn{e_4}{w}^2+\inn{e_3}{w}^2$, this inequality holds
trivially. $\Box$ 

\bigskip

We now modify Shi-Tam's proof of Theorem \ref{thm:shi-tam}.
Let $\Si^\al, \Si^\al_r, E^\al,g^\al=h^2 dr^2+ g_r$
be defined as in in Section \ref{sec:shi-tam}.
Shi and Tam's proof in \cite{Sh-Ta} shows that there is a 
unique solution to (\ref{eqn:shi-tam-foliation})
on $E^\al\cong \Si^\al\times[0,\infty)$
with the initial condition
\begin{equation}\label{initial}
h(x,0)=\frac{H_0^\al}{\bH-\inn{X}{\bnu} }
\end{equation}
such that
\[
|h(x,r)-1|\leq \frac{C}{r}
\]
for $r\geq 1$. Equip $E^\alpha$ with the metric
$g^\al = h^2 dr^2 + g_r$.
Then $g^\al$ has zero scalar curvature, and
the mean curvature of $\Si^\al$ in $(E^\al,g^\al)$ is
$\bH-\inn{X}{\bnu}$. 

Define 
$$
m^\al(r)=\frac{1}{8\pi G}\int_{\Si^\al_r}(H_0-H)d\si_r,\ \ 
m_\infty^\al=\lim_{r\to\infty} m^\al(r)
$$
as in Section \ref{sec:shi-tam}. By monotonicity (\ref{eqn:monotone}), 
$m_\infty^\al\leq m^\al(0)$. By Lemma \ref{thm:XHP}, 
$$
m^\al(0)=\frac{1}{8\pi G}\int_\Si(H^\al_0-(\bH-\inn{X}{\bnu}))d\si
      \leq\frac{1}{8\pi G}\int_\Si(H^\al_0-\sqrt{H^2-P^2})d\si = E(\Si^\al).
$$

\subsection{Asymptotically flat Lipschitz metric}\label{sec:metric}

Following \cite{Sh-Ta}, we glue $(E^\al,g^\al)$ in Section \ref{sec:end} to 
$(\Om,\hg)$ in Section \ref{sec:interior} along $\Si^\al$ to obtain a complete 
noncompact three manifold $M$ with a continuous Riemannian metric $\tg$ such that
\begin{enumerate}
\item $\tg$ is smooth on $M\setminus\Om$ and $\bOm$, and is Lipschitz near $\pa\Om$.
\item Each end $E^\al$ of $M$ is asymptotically Euclidean.
\item The scalar curvature $R$ of $M\setminus\pa\Om$ is nonnegative and
      is in $L^1(M)$.
\end{enumerate}
The mean curvatures of $\Si$ with respect to $\hg=\tg|_\Om$
is $\hH=\bH+4\bnu(u)$, and the mean curvature
of $\Si$ with respect to $g^\al=\tg|_{E^\al}$ is $\bH-\inn{X}{\bnu}$. 

Let $\tnu$ be the outward unit normal of $\pa\Om$ with respect to $\tg$.
There exists $\ep>0$ such that 
$(x,t) \mapsto \exp_x(t\tilde{\nu}(x))$
defines an open embedding $i:\pa\Om\times(-\ep,\ep)\to M$.
The image $T=i(\pa\Om\times (-\ep,\ep))$ is a tubular neighborhood of $\pa\Om$ in $M$.
We use the smooth structure on $\pa\Om\times (-\ep,\ep)$ to define
the smooth structure on $T$. We have
$$
i^*\tg = dt^2 + \rho_{ij}(x,t)dx_i dx_j
$$
where $(x_1,x_2)$ are local coordinates on $\pa\Om$.

We choose local orthonormal frame $e_1,e_2$ of $\pa\Om$ and
parallel transport them along $\frac{\pa}{\pa t}$. Then
$e_1,e_2,$ and $e_3=\frac{\pa}{\pa t}$ form a local orthonormal frame
on $\left(\pa\Om\times(-\ep,\ep),i^*\tg\right)$ such that
\begin{enumerate}
\item $e_1,e_2$ are tangent to slides $\Si_t=\pa\Om\times \{ t\}$. 
\item $e_3$ is normal to $\Si_t$.
\item $\nabla_3 e_i=0$, where $\nabla$ is the Levi-Civita
     connection with respect to $i^*\tg$.
\end{enumerate}
The mean curvature of $\Si_t$ in $\pa\Omega \times (-\ep,\ep)$ defines
a function $H$ on $\pa\Omega\times (-\ep,\ep)$ which is discontinuous
at $t=0$ and smooth away from $t=0$. 

Note that $(M,\tg)$ is uniquely determined by $(\Om,\hg)$, which 
is uniquely determined by $(\Om,\bg)$.
As explained in Section \ref{sec:general}, our main result
Theorem \ref{thm:positive} follows from the following positive
mass theorem of $(M,\tg)$.
\begin{thm}\label{thm:tgmass}
The ADM mass $m_\infty^\al$ of the end $(E^\al,g^\al)$ is nonnegative for
$\al=1,\ldots,\ell$, and $m^\al_\infty=0$ for some $\al$ if and only if
$\ell=1$ and $M$ is the Euclidean space.
\end{thm}

We will prove Theorem \ref{thm:tgmass} in Section \ref{sec:tgmass}.

\subsection{Dirac spinor} \label{sec:dirac}
In the rest of this paper,
$L^p, L^p_\loc, L^p_0, W^{k,p}, W^{k,p}_\loc, W^{k,p}_0$ are defined
as in \cite[Chapter 7]{Gi-Tr}.

The spinor bundle $S$ over $M$ is a trivial complex vector bundle of
rank $2$. Let $\vD:C^\infty(M,S)\to C^\infty(M,S)$ be the Dirac operator defined by 
the Levi-Civita connection of $\tg$. It can be extended to 
$\vD:W^{1,2}_\loc(M,S)\to L^2_\loc(M,S)$.  

Let $T$ be the tubular neighborhood of $\pa\Om$ in $M$ defined in Section \ref{sec:metric}.
We identify $T$ with $\pa\Om\times (-\ep,\ep)$ and study
the Dirac operator on $\pa\Om\times (-\ep,\ep)$. 
Let $e_1,e_2,e_3$ be defined as in Section \ref{sec:metric}, and
let $\theta^1,\theta^2,\theta^3$ be the dual coframe. Then $\theta^3=dt$. Let
$$
\vD^t= c(\theta^3)\left(c(\theta^1)\na^t_1 + c(\theta^2) \na^t_2\right)
$$
where $\na^t$ is the Levi-Civita connection on $\Si_t$. Then
$$
\vD \psi =c(\theta^3)\left(\frac{\pa\psi}{\pa t} - \vD^t \psi
+\frac{H}{2}\psi\right).
$$

Let
$$
\vD'= \vD-\frac{1}{2}\beta(t)H c(\theta^3)
$$
where $\beta:(-\ep,\ep)\to \bR$ is a smooth function such that $\beta(t)=\beta(-t)$
and
$$
\beta(t)=\left\{\begin{array}{ll} 
1 & |t| \leq \ep/3\\
0 & |t|> 2\ep/3\end{array}
\right.
$$
Then $\vD'$ extends to a first-order differential operator on $M$
with smooth coefficients.

In Section \ref{sec:exist}, we will prove the following existence and uniqueness of Dirac
spinor with prescribed asymptotics. 
\begin{thm}\label{thm:dirac}
Let $\psi^1,\ldots,\psi^\ell$ be constant spinors defined on 
the ends $E^1,\ldots,E^\ell$. Then there exists a unique 
spinor $\psi\in W^{1,2}_\loc(M,S)$ such that
\begin{enumerate}
\item $\vD\psi=0$.
\item $\psi\in C^\infty(M\setminus \pa\Om, S)$.
\item $\psi\in W^{1,p}_\loc(M,S)$ for any $2\leq p <\infty$.
\item On each end $E^\al$, let $\rho$ be defined as in
  Theorem \ref{thm:unique}. Then
$$
\lim_{\rho\to\infty}\rho^{1-\ep}|\psi-\psi^\al|=0
$$
for any $\ep>0$.
\end{enumerate} 
\end{thm}

\medskip

In general, the mean curvature
along $\pa\Om$ is discontinuous, so the Dirac spinor $\psi$ in 
Theorem \ref{thm:dirac} is not in $C^1(M,S)$.  

\subsection{Boundary values of $W^{1,2}$ functions}
\label{sec:boundary}
We recall the following result from \cite{Mo}, where
$H_p^m$ corresponds to $W^{m,p}$ in our notation.
\begin{thm}[{\cite[Theorem 3.4.5]{Mo}}]\label{thm:boundary}
If $G$ is bounded and of class $C^{m-1}_1$ the functions 
$u\in C^{m-1}_p(\bar{G})$ are dense in any space $H^m_p(G)$ with $p\geq 1$
and there is a bounded operator $B$ from $H_p^m(G)$ into $H_p^{m-1}(\pa G)$
such that $Bu=u|_{\pa G}$ whenever $u\in C_1^{m-1}(\bar{G})$. If $u_n \to u$
in $H_p^m(G)$, then $u_n\to u$ in $H_p^{m-1}(\pa G)$. If $p>1$, the mapping is
compact.
\end{thm}

Recall that we have translate $\Si^\al_0\subset \bR^3$ such that
there is $a>0$ such that the closed ball $B_a$ of radius $a$ centered at
the origin is disjoint from $E^\al$. Choose $L>0$ such that
$E^\al$ contains $\bR^3\setminus B_L$. For $r>L$, let
$$
E^\al_r= E^\al\setminus B_r,\ \ 
S^\al_r=\pa E^\al_r,\ \ 
M_r= M\setminus\cup_{\al=1}^\ell E^\al_r.
$$
For a fixed $r>L$, let $G_+$ be the interior of $\Omega$, and let
$G_-$ be the interior of $M_r\setminus \Omega$.
We have a disjoint union
$$
M_r= G_+ \cup \pa\Om \cup G_- \cup \pa M_r,
$$
where
$$
\bar{G}_+= G_+\cup \pa\Om,\ \
\bar{G}_-= \pa\Om \cup G_- \cup \pa M_r.
$$ 
Let $r_\pm:W^{1,2}(M)\to W^{1,2}(G_\pm)$ be the restriction map, and
let $b_\pm:W^{1,2}(G_\pm) \to L^2(\pa\Om)$ be the bounded
linear operator in Theorem \ref{thm:boundary}. Let
$B_\pm=b_\pm\circ r_\pm: W^{1,2}(M) \to L^2(\pa\Om)$.
Given $u\in W^{1,2}(M)$, there exists a sequence 
$\{u_n\}\subset C^{m-1}_1(M_r)$ such that $u_n\to u|_{\Om\cup G_-}$
in $W^{1,2}(\Om\cup G_-)$.
Then $B_+ u = B_- u =\lim_{n\to \infty}(u_n|_{\pa\Om})$.

\subsection{Estimates near $\Omega$}\label{sec:near}
Let $L>0$ be chosen as in Section \ref{sec:boundary}. For $r>L$, 
let $M_r, S_r^\al$ be  defined as in Section \ref{sec:boundary}.
The goal of this subsection is to establish the following estimate
which will be a crucial ingredient in the proof of Theorem
\ref{thm:dirac}. 
\begin{pro}\label{thm:key}
For $r>L$ and $\psi\in W^{1,2}_\loc(M,S)\cap C^\infty(M\setminus M_L,S)$, we have
\begin{equation}\label{eqn:key}
2\int_{M_r}|\vD\psi|^2 \geq
\frac{1}{10}\int_{M_r}|\na\psi|^2+\frac{1}{16}\int_\Om u^{-2}|du|^2|\psi|^2
+\sum_{\al=1}^\ell\int_{S^\al_r}\inn{\frac{H}{2}\psi-c(\nu)\vvD\psi}{\psi}
\end{equation}
where $\vvD$ is the Dirac operator on $S^\al_r$.
\end{pro}

\begin{lm}\label{thm:Usurface}
Let $U$ be an open set of $M$. For any spinor $\eta\in W^{1,2}_0(U,S)$, 
$\psi\in W^{1,2}_\loc(U,S)$, we have
\begin{equation}\label{eqn:Uadj}
\int_U\inn{\vD\psi}{\eta}=\int_U \inn{\psi}{\vD\eta}
\end{equation}
\begin{equation}\label{eqn:Uwei}
\int_U\inn{\vD\psi}{\vD\eta}
=\int_U\inn{\na\psi}{\na\eta}_\tg
+\int_{\pa\Om\cap U}(2\bnu(u)+\frac{1}{2}\inn{X}{\bnu})\inn{\psi}{\eta}
\end{equation}
\end{lm}
\paragraph{Proof.} 
By the discussion in Section \ref{sec:boundary},
the right hand side of (\ref{eqn:Uwei}) makes sense.
It suffices to show that (\ref{eqn:Uadj}) and (\ref{eqn:Uwei})
hold for $\psi\in C^\infty(U,S),\eta\in C^\infty_0(U,S)$.

Let $U_1=\Om\cap U$, $U_2= U\setminus\bOm$, and
$I=\pa\Omega\cap U$. 
We have 
\begin{eqnarray}
\int_{U_1}\inn{\vD\psi}{\eta}&=&
\int_{U_1}\inn{\psi}{\vD\eta}
+\int_I\inn{c(\hnu)\psi}{\eta}
\label{eqn:aI}\\
\int_{U_2}\inn{\vD\psi}{\eta}&=&
\int_{U_2}\inn{\psi}{\vD\eta}
+\int_I\inn{-c(\hnu)\psi}{\eta}
\label{eqn:aII}
\end{eqnarray}
where $\hnu$ is the outward unit normal of $\pa\Om$ in $(\Om,\hg)$
(see e.g. \cite[Proposition 3.4]{Bo-Wo}). Equation (\ref{eqn:Uadj}) is
the sum of (\ref{eqn:aI}) and (\ref{eqn:aII}).

We also have 
\begin{eqnarray}
\int_{U_1}\inn{\vD\psi}{\vD\eta}&=&
\int_{U_1}\inn{\na\psi}{\na\eta}_\tg
+\int_I\inn{\frac{1}{2}(\bH+4\bnu(u))\psi-c(\hnu)\vvD\psi}{\eta}
\label{eqn:wI}\\
\int_{U_2}\inn{\vD\psi}{\vD\eta}&=&
\int_{U_2}\inn{\na\psi}{\na\eta}_\tg
+\int_I\inn{\frac{1}{2}(\inn{X}{\bnu}-\bH)\psi+c(\hnu)\vvD\psi}{\eta}
\label{eqn:wII}
\end{eqnarray}
where $\vvD$ is the  Dirac operator on $\pa\Om$ (see
e.g. \cite{HMZ}). Equation (\ref{eqn:Uwei}) is the sum of 
(\ref{eqn:wI}) and (\ref{eqn:wII}). $\ \ \Box$

The proof of Lemma \ref{thm:Usurface} also gives the following:
\begin{lm}\label{thm:Msurface}
For $r>L$ and $\psi,\eta\in W^{1,2}_\loc(M,S)$, we have
\begin{equation}\label{eqn:Madj}
\int_{M_r}\inn{\vD\psi}{\eta}=\int_{M_r} \inn{\psi}{\vD\eta}
+\sum_{\al=1}^\ell\int_{S^\al_r}\inn{c(\nu)\psi}{\eta}
\end{equation}
For $r>L$, $\psi\in W^{1,2}_\loc(M,S)\cap C^\infty(M\setminus
M_L,S)$, and $\eta\in W^{1,2}_\loc(M,S)$, we have
\begin{equation}\label{eqn:Mwei}
\int_{M_r}\inn{\vD\psi}{\vD\eta}
=\int_{M_r}\inn{\na\psi}{\na\eta}_\tg
+\int_{\pa\Om}(2\bnu(u)+\frac{1}{2}\inn{X}{\bnu})\inn{\psi}{\eta}
+\sum_{\al=1}^\ell\int_{S^\al_r}\inn{\frac{H}{2}\psi-c(\nu)\vvD\psi}{\eta}
\end{equation}
\end{lm}

\begin{lm}\label{thm:DvD}
Let $\bar{\vD}$ and $\vD$ denote the Dirac operators on $S|_\Omega$ defined
by the Levi-Civita connections of $\bg$ and $\hg=u^4\bg$,
respectively. Then
\begin{equation}\label{eqn:DvD}
\vD\psi = \frac{1}{u^2}\bar{\vD}\psi + \frac{1}{2u^3}\bc(du)\psi
\end{equation}
where $\bc$ is the Clifford multiplication defined by $\bg$.
\end{lm}
\paragraph{Proof}
The tangent bundle of $\Omega$ is trivial, so there
exists global orthonormal frame $\{\be_1,\be_2,\be_3\}$
with respect to $\bg$. Let $\{\bth^1,\bth^2,\bth^2\}$
be the dual coframe. We have
\begin{eqnarray*}
\bna_{\be_i} \be_j &=& \bGa_{ij}^k\be_k\\
d\bth^i&=&-\bGa_{kj}^i\bth^k\wedge\bth^j\\
\bna_{\be_i}\psi&=&\be_i(\psi)+\frac{1}{4}\bGa_{ij}^k
\bc(\bth^j)\bc(\bth^k)\psi
\end{eqnarray*}

Let 
$$
e_i=u^{-2} \be_i,\ \  \theta^i=u^2 \bth^i.
$$ 
Then $\{e_1,e_2,e_3\}$ is a global orthonormal frame
with respect to $\hg$, and $\{\theta^1,\theta^2,\theta^3\}$
is the coframe. Let $u_i=\be_i(u)$, $\psi_i=\bna_{\be_i}\psi$. 
Then
\begin{eqnarray*}
d\theta^i&=&  2u_j u\bth^j\wedge \bth^i 
-\sum_{j,k}u^2\bGa_{kj}^i\bth^k\wedge\bth^j
=-\left(\bGa^i_{kj}\bth^k+\frac{2u_j}{u}\bth^i\right)\wedge \theta^j\\
\Ga^i_{kj}&=&
\left(\bGa_{lj}^i\bth^l + \frac{2u_j}{u}\bth^i\right)(e_k)
=u^{-2}\bGa^i_{kj}+ 2u^{-3}u_j\delta_{ik}\\
\na_{e_i}\psi&=&
e_i(\psi)+\frac{1}{4}\Ga_{ij}^k c(\theta^j)c(\theta^k)\psi
=\frac{1}{u^2}\psi_i + \frac{1}{2u^3}\bc(du)\bc(\bth^i)\psi
\end{eqnarray*}
Note that $\bc(\bth^i)=c(\theta^i)$, so
$$
\vD\psi =c(\theta^i)\na_{e_i}\psi
        =\frac{1}{u^2}\bc(\bth^i)\psi_i +
       \frac{1}{2u^3}\bc(\bth^i)\bc(du)\bc(\bth^i)\psi
       =\frac{1}{u^2}\D\psi +
       \frac{1}{2u^3}\bc(du)\psi.\ \ \ \Box
$$

\begin{lm}\label{thm:dpsi}
\begin{equation}\label{eqn:hgbg}
|\na\psi|^2_\hg = \frac{1}{u^4}|\bna\psi|^2_\bg
+\frac{3}{4u^6}|du|_\bg^2|\psi|^2
+\frac{1}{u^5}\mathrm{Re}\inn{\D \psi}{\bc(du)\psi}
-\frac{2}{u^5}\mathrm{Re}\inn{u_i\psi_i}{\psi}
\end{equation}
\begin{equation}\label{eqn:bghg}
|\bna\psi|^2_\bg= u^4|\na\psi|^2_\hg
+\frac{3u^2}{4}|du|_\hg^2|\psi|^2
-u^3\mathrm{Re}\inn{\vD\psi}{c(du)\psi}
+2u^3\mathrm{Re}\inn{e_i(u)\na_{e_i}\psi}{\psi}
\end{equation}
\end{lm}

\paragraph{Proof}
We use the notation in the proof of Lemma \ref{thm:DvD}, 
where we calculated that 
\begin{equation}\label{eqn:nabna}
\na_{e_i}\psi=\frac{1}{u^2}\psi_i +\frac{1}{2u^3}\bc(du)\bc(\bth^i)\psi 
\end{equation}
Note that $\inn{\psi_1}{\psi_2}$ does not depend on the metric on 
the tangent bundle of $\Omega$. We have
\begin{equation}\label{eqn:square}
|\na_{e_i}\psi|^2 =\frac{1}{u^4}|\psi_i|^2 
+ \frac{|du|^2_\bg}{4u^4}|\psi|^2
+\frac{1}{u^5}\mathrm{Re}\inn{\psi_i}{\bc(du)\bc(\bth^i)\psi}.
\end{equation}
Note that the last term on the right hand side of (\ref{eqn:square})
is {\em not} a sum, and can be rewritten as follows.
\begin{eqnarray*}
\inn{\psi_i}{\bc(du)\bc(\bth^i)\psi}
&=& \inn{\psi_i}{u_j\bc(\bth^j)\bc(\bth^i)\psi}\\
&=& \inn{\psi_i}{u_j(-\bc(\bth^i)\bc(\bth^j)-2\delta_{ij})\psi}\\
&=& \inn{\bc(\bth^i)\psi_i}{u_j\bc(\bth^j)\psi}
    -2u_i\inn{\psi_i}{\psi}
\end{eqnarray*}

We now sum over $i=1,2,3$ and obtain
\begin{eqnarray*}
|\na\psi|^2_\hg&=& \sum_{i=1}^3 |\na_{e_i}\psi|^2\\
&=& \frac{1}{u^4}|\bna\psi|_\bg^2+\frac{3|du|^2_\bg}{4u^6}|\psi|^2
+\frac{1}{u^5}\mathrm{Re}\left(\inn{\D\psi}{\bc(du)\psi}
-2u_i\inn{\psi_i}{\psi}\right)
\end{eqnarray*}
This proves (\ref{eqn:hgbg}). By symmetry, we have
$$
|\bna\psi|^2_\bg = u^4|\na\psi|^2_\hg
+\frac{3u^6}{4}|d(u^{-1})|^2_g|\psi|^2
+u^5\mathrm{Re}\inn{\vD\psi}{c(d(u^{-1}))\psi}
-2u^5\mathrm{Re}\inn{e_i(u^{-1})\na_{e_i}\psi}{\psi}
$$
which is equivalent to (\ref{eqn:bghg}). $\ \ \Box$

\begin{lm}\label{thm:onetwothree} 
Let $d\bsi$ denote the volume form of $\bg$, and
$\mathrm{div}$ denote the divergence defined by $\bg$. Then
\begin{equation}\label{eqn:one}
\int_\Om |\na\psi|^2_\hg d\si
=\int_\Om \left(u^2|\bna\psi|^2_\bg 
+\frac{3}{4}|du|^2_\bg|\psi|^2
+u\mathrm{Re}\inn{\D\psi}{\bc(du)\psi} 
-2u\mathrm{Re}\inn{u_i\psi_i}{\psi}\right)d\bsi
\end{equation}
\begin{equation}\label{eqn:two}
2\int_{\pa\Om}\bnu(u)|\psi|^2
\geq \int_\Om \left(2|du|^2_\bg |\psi|^2
+\frac{1}{2}|X|^2_\bg u^2|\psi|^2
-\frac{1}{2}(\mathrm{div}X) u^2|\psi|^2 
 +4 u\mathrm{Re}\inn{u_i\psi_i}{\psi}\right) d\bsi
\end{equation}
\begin{equation}\label{eqn:three}
\frac{1}{2}\int_{\pa\Om}\inn{X}{\bnu}|\psi|^2
\geq \int_\Om\left(\frac{1}{2}(\mathrm{div}X)u^2|\psi|^2
-\frac{1}{2}|X|_\bg^2u^2|\psi|^2
-\frac{|\bna(u^2|\psi|^2)|^2_\bg}{8u^2|\psi|^2}\right)d\bsi
\end{equation}
\end{lm}

\paragraph{Proof} We have $d\bsi=u^6 d\si$, so (\ref{eqn:one}) follows from (\ref{eqn:hgbg}).
To prove (\ref{eqn:two}), note that
\begin{eqnarray*}
2\int_{\pa\Om}\bnu(u)|\psi|^2 d\si
&=& 2\int_\Om \mathrm{div}(u\bna u|\psi|^2)d\bsi\\
&=& 2\int_\Om|\bna u|_\bg^2|\psi|^2 d\bsi
    + 2\int_\Om u\bar{\Delta} u |\psi|^2 d\bsi
   + 2\int_\Om u\inn{\bna u}{\bna (|\psi|^2)}_\bg d\bsi
\end{eqnarray*}
where
$$
\bar{\Delta}u=\frac{\bar{R}}{8}u\geq 
\frac{1}{4}(|X|^2_\bg-\mathrm{div}X)u,
$$
$$
\inn{\bna u}{\bna (|\psi|^2)}_\bg = u_i \be_i(|\psi|^2)
=2\mathrm{Re}\inn{u_i\psi_i}{\psi},
$$
so (\ref{eqn:two}) holds. Finally,
\begin{eqnarray*}
\frac{1}{2}\int_{\pa\Om}\inn{X}{\bnu}|\psi|^2
&=&\frac{1}{2}\int_\Om \mathrm{div}(X u^2 |\psi|^2) d\bsi\\
&=&\frac{1}{2}\int_\Om (\mathrm{div}X) u^2 |\psi|^2 d\bsi
  +\frac{1}{2}\int_\Om\inn{X}{\bna(u^2|\psi|^2)}_\bg d\bsi\\
&=&\frac{1}{2}\int_\Om (\mathrm{div}X) u^2 |\psi|^2 d\bsi
  +2\int_\Om\inn{ \frac{1}{\sqrt{2}}Xu|\psi| }{
    \frac{\bna(u^2|\psi|^2_)}{\sqrt{8} u|\psi|} }_\bg d\bsi\\
&\geq&\frac{1}{2}\int_\Om (\mathrm{div}X) u^2 |\psi|^2 d\bsi
-\frac{1}{2}\int_\Om |X|^2_\bg u^2|\psi|^2 d\bsi
-\frac{1}{8}\int_\Om\frac{|\bna(u^2|\psi|^2)|^2_\bg}{u^2|\psi|^2}d\bsi
\end{eqnarray*}
This proves (\ref{eqn:three}). $\Box$

\begin{lm}\label{thm:onetenth}
$$
\int_\Om|\na\psi|^2_\hg d\si+2\int_{\pa\Om}\bnu(u)|\psi|^2
+\frac{1}{2}\int_{\pa\Om}\inn{X}{\bnu}|\psi|^2\\
\geq\int_\Om \left(\frac{1}{10}|\na\psi|^2_\hg 
+ \frac{1}{16u^2}|du|^2_\hg|\psi|^2-|\vD\psi|^2\right)d\si 
$$
\end{lm}
\paragraph{Proof}
By Lemma \ref{thm:onetwothree}, we have
\begin{eqnarray*}
&&\int_\Om|\na\psi|^2_\hg d\si+2\int_{\pa\Om}\bnu(u)|\psi|^2
+\frac{1}{2}\int_{\pa\Om}\inn{X}{\bnu}|\psi|^2\\
&\geq& \int_\Om \left(u^2 |\bna\psi|_\bg^2
+\frac{11}{4}|du|_\bg^2|\psi|^2 
+2 u\mathrm{Re}\inn{u_i\psi_i}{\psi}
-\frac{1}{8}\frac{|\bna(u^2|\psi|^2)|^2 }{u^2|\psi|^2 } 
+ u\mathrm{Re}\inn{\vD\psi}{\bc(du)\psi}\right)d\bsi
\end{eqnarray*}

We rewrite
$$
-\frac{1}{8}\frac{|\bna(u^2|\psi|^2)|^2 }{u^2|\psi|^2 }
$$
as follows:
\begin{eqnarray*}
\bna(u^2|\psi|^2)&=&2u\bna u|\psi|^2 + u^2
2\mathrm{Re}\inn{\bna\psi}{\psi}\\
|\bna(u^2|\psi|^2)|^2 &=& 4u^2|du|^2_\bg|\psi|^4 + 
4u^4\sum_{i=1}^3(\mathrm{Re}\inn{\psi_i}{\psi})^2
+8u^3|\psi|^2 \mathrm{Re}\inn{u_i\psi_i}{\psi}\\
-\frac{|\bna(u^2|\psi|^2)|^2}{8 u^2|\psi|^2}
&=&-\frac{1}{2}|d u|_\bg^2|\psi|^2
   -\frac{1}{2}\frac{u^2}{|\psi|^2}\sum_{i=1}^3(\mathrm{Re}\inn{\psi_i}{\psi})^2
   -u\mathrm{Re}\inn{u_i\psi_i}{\psi}\\
&\geq &-\frac{1}{2}|d u|_\bg^2|\psi|^2
   -\frac{1}{2}u^2|\bna \psi|_\bg^2
   -u\mathrm{Re}\inn{u_i\psi_i}{\psi}
\end{eqnarray*}

So
\begin{equation}\label{eqn:sumthree}
\int_\Om|\na\psi|^2_\hg d\si+2\int_{\pa\Om}\bnu(u)|\psi|^2
+\frac{1}{2}\int_{\pa\Om}\inn{X}{\bnu}|\psi|^2 \hspace{1in}
\end{equation}
$$
\geq \int_\Om \left(\frac{1}{2}u^2 |\bna\psi|_\bg^2
+\frac{9}{4}|du|_\bg^2|\psi|^2 
+u\mathrm{Re}\inn{u_i\psi_i}{\psi} 
+ u\mathrm{Re}\inn{\D\psi}{\bc(du)\psi} \right)d\bsi.
$$

In the rest of this proof, we will write 
$$
|du|^2=|du|_\hg^2,\ \ 
|\na\psi|^2=|\na\psi|^2_\hg.
$$
By (\ref{eqn:bghg}),
\begin{equation}\label{eqn:I}
\frac{1}{2}u^2|\bna\psi|^2_\bg
= \frac{u^6}{2}|\na\psi|^2
 +\frac{3u^4}{8}|du|^2|\psi|^2
 -\frac{u^5}{2}\mathrm{Re}\inn{\vD\psi}{c(du)\psi}
+u^5\mathrm{Re}\inn{e_i(u)\na_{e_i}\psi}{\psi}.
\end{equation}
By (\ref{eqn:nabna}) and symmetry,
\begin{equation}\label{eqn:bnana}
\psi_i=u^2\na_{e_i}\psi+\frac{u^3}{2}c(d(u^{-1}))c(\theta^i)\psi
=u^2\na_{e_i}\psi-\frac{u}{2}c(du)c(\theta^i)\psi.
\end{equation}
So
\begin{eqnarray*}
u\mathrm{Re}\inn{u_i\psi_i}{\psi}
&=&u\mathrm{Re}
\inn{u^2e_i(u)\left(u^2\na_{e_i}\psi-\frac{u}{2}c(du)c(\theta^i)\right)\psi}{\psi}\\
&=&u^5\mathrm{Re}\inn{e_i(u)\na_{e_i}\psi}{\psi}-\frac{u^4}{2}\mathrm{Re}\inn{c(du)^2\psi}{\psi}
\end{eqnarray*}
which implies
\begin{equation}\label{eqn:III}
u\mathrm{Re}\inn{u_i\psi_i}{\psi}
=u^5\mathrm{Re}\inn{e_i(u)\na_{e_i}\psi}{\psi}+\frac{u^4}{2}|du|^2|\psi|^2
\end{equation}
We also have
\begin{eqnarray*}
u\mathrm{Re}\inn{\D\psi}{\bc(du)\psi}
&=& u\mathrm{Re}\inn{u^2 \vD\psi-\frac{1}{2u}\bc(du)\psi}{\bc(du)\psi}\\
&=& u^5\mathrm{Re}\inn{\vD\psi}{c(du)\psi}-\frac{1}{2} u^4\inn{c(du)\psi}{c(du)\psi}
\end{eqnarray*}
which implies
\begin{equation}\label{eqn:IV}
u\mathrm{Re}\inn{\D\psi}{\bc(du)\psi}
=u^5\mathrm{Re}\inn{\vD\psi}{c(du)\psi}-\frac{1}{2}u^4|du|^2|\psi|^2.
\end{equation}
Let $v=\log(u)$. Then
(\ref{eqn:sumthree}), (\ref{eqn:I}), (\ref{eqn:III}), (\ref{eqn:IV})
imply that
\begin{eqnarray*}
&&\int_\Om|\na\psi|^2_\hg d\si+2\int_{\pa\Om}\bnu(u)|\psi|^2
+\frac{1}{2}\int_{\pa\Om}\inn{X}{\bnu}|\psi|^2\\ 
&\geq&\int_\Omega\left(\frac{1}{2}|\na\psi|^2
+\frac{21}{8}|dv|^2|\psi|^2
+\frac{1}{2}\mathrm{Re}\inn{\vD\psi}{c(dv)\psi}
+2\mathrm{Re}\inn{e_i(v)\na_{e_i}\psi}{\psi}\right)d\si\\
&\geq&\int_\Omega\left(\frac{1}{2}|\na\psi|^2
+\frac{21}{8}|dv|^2|\psi|^2
-(|\vD\psi|^2+\frac{1}{16}|dv|^2|\psi|^2)
-(\frac{2}{5}|\na\psi|^2 +\frac{5}{2}|dv|^2|\psi|^2\right)d\si\\
&\geq&\int_\Omega\left(\frac{1}{10}|\na\psi|^2
+\frac{1}{16}|dv|^2|\psi|^2-|\vD\psi|^2\right)d\si\ \ \Box
\end{eqnarray*}

\medskip

\paragraph{\bf Proof of Proposition \ref{thm:key}}
By (\ref{eqn:Mwei}),
\begin{eqnarray*}
\int_{M_r}|\vD\psi|^2& =& \int_{M_r\setminus\Om}|\na\psi|^2 +\int_\Om|\na\psi|^2 
+ \int_{\pa\Om}(2\bnu(u)+\frac{1}{2}\inn{X}{\bnu})|\psi|^2\\
&&+\sum_{\al=1}^\ell \int_{S^\al_r}\inn{\frac{H}{2}\psi-c(\nu)\vvD\psi}{\psi},
\end{eqnarray*}
where
$$
 \int_\Om|\na\psi|^2
+\int_{\pa\Om}(2\bnu(u)+\frac{1}{2}\inn{X}{\bnu})|\psi|^2
\geq \frac{1}{10}\int_\Om |\na\psi|^2 +\frac{1}{16}\int_\Om
u^{-2}|du|^2|\psi|^2 -\int_\Om|\vD\psi|^2
$$
by Lemma \ref{thm:onetenth}. Therefore,
\begin{eqnarray*}
&& \int_{M_r}|\vD\psi|^2 +\int_{\Om}|\vD\psi|^2\\
 &\geq &\int_{M_r\setminus\Om}|\na\psi|^2
+\frac{1}{10}\int_\Om|\na\psi|^2 
+\frac{1}{16}\int_\Om u^{-2}|du|^2|\psi|^2 
+  \sum_{\al=1}^\ell \int_{S^\al_r}\inn{\frac{H}{2}\psi-c(\nu)\vvD\psi}{\psi}
\end{eqnarray*}
which implies (\ref{eqn:key}).$\ \ \Box$

\subsection{Proof of Theorem \ref{thm:dirac}}\label{sec:exist}
We modify Parker and Taubes's proof of \cite[Theorem 4.1]{Pa-Ta}.
In \cite{Pa-Ta}, the hypersurface Dirac operator is studied,
so the second fundamental form of the spacelike hypersurface in the
spacetime is involved in estimates. Here we consider the Riemannian
case, so the estimates are simpler in certain steps. The main
difficulty in our case comes from  the discontinuity of the zeroth 
order term of the Dirac operator along $\pa\Om$.

Let $L>0$ be chosen as in Section \ref{sec:boundary}. For $r>L$, let
$E^\al_r, M_r, S^\al_r$ be defined as in Section \ref{sec:boundary}.
Choose a smooth function $\beta^\al_L$ on each end $E^\al$ such that 
(i) $0\leq \beta^\al_L\leq 1$, (ii) $\beta^\al_L\equiv 1$ on $E^\al_{3L}$,
(iii) $\beta^\al_L\equiv 0$ on $E^\al\cap M_{2L}$, and 
(iv) $|\na\beta^\al_L|\leq 2/L$. 
Then $\beta_L^\al \psi^\al$ extends to a smooth section of $S$ over
$M$. Define
\begin{equation}\label{eqn:ends}
\psi_0=\sum_{\al=1}^\ell \beta_L^\al \psi^\al\in C^\infty(M,S). 
\end{equation}
We wish to find $\psi_1\in W^{1,2}(M,S)$ such that 
\begin{equation}
\vD\psi_1=-\vD\psi_0
\end{equation}
and
$$
\lim_{r\to\infty}r^{1-\ep}|\psi_1|=0.
$$
Then $\psi=\psi_0+\psi_1$ is the desired solution.

Let $U$ be an open subset of $M$, and let $\psi\in L^2_\loc(U,S)$.
Given $\eta\in C^\infty(U,S)$, $\psi$ is said to satisfy
$$
\vD\psi=\eta
$$
in the weak sense if
$$
\int_U\inn{\psi}{\vD\phi}=\int_U \inn{\eta}{\phi}
$$
for any $\phi\in C^1_0(U,S)$, or equivalently, for any
$\phi\in W^{1,2}_\loc(U,S)$.

\begin{lm}\label{thm:regularity}
Let $U$ be an open subset of $M$, and let $\psi\in L^2_\loc(U,S)$ be a
weak solution to
$$
\vD\psi=\eta
$$
where $\eta\in C^\infty(U,S)$. Then $\psi\in
C^\infty(U\setminus\pa\Om,S)$, and $\psi\in W^{1,p}_\loc(U,S)$
for any $2\leq p<\infty$.
\end{lm}

\paragraph{Proof} Recall from Section \ref{sec:dirac} that $\vD=\vD'+A$, where $\vD'$ is a first order
elliptic operator with smooth coefficients, and 
$A\in L^\infty(U,\mathrm{End}(S))\cap C^\infty(U\setminus \pa\Om,\mathrm{End}(S))$.  
So $\psi\in C^\infty(U\setminus\pa\Om,S)$ by elliptic regularity. Now
$$
\vD'\psi=A\psi+ \eta
$$
where $A\psi+\eta\in L^2_\loc(U,S)$. By elliptic regularity, 
$\psi\in W^{1,2}_\loc(U)\subset L^6_\loc(U)$, so $A\psi+\eta \in
L^6_\loc(U)$. By elliptic regularity again, $\psi\in
W^{1,6}_\loc(U,S)\subset C^0(U,S)$. 
So $A\psi+\eta\in L^p_\loc(U,S)$ for any $p\geq 2$, which implies
$\psi\in W^{1,p}_\loc(U,S)$ for any $p\geq 2$. $\ \ \Box$

\medskip

We recall some weighted  Sobolev spaces introduced
in \cite{Pa-Ta}. The distance function from the 
origin is a smooth function $\rho: E^\al_L\to \bR$
such that $\rho(E^\al_L)=[L,\infty)$.This defined a smooth function $\rho$ on $M_L$
such that $\rho^{-1}(r)=\pa M_r$ for $r\geq L$.
Fix a smooth function $\si$ on $M$ such that
(i) $\si\geq 1$, (ii) $\si=\rho$ on $E^\al_{2L}$ and  
(iii) $\si=1$ in $M_L$.
\begin{df}
Given a pair $(\dep)$ such that
$$
p\geq 2,\ \
\frac{1}{2}-\frac{3}{p}\leq \delta \leq 2-\frac{3}{p},
$$
let $W^{0,p}_\delta$, $W^{1,p}_\delta$ be the completions of $C_0^\infty(M,S)$
with respect to the norms
\begin{equation}
\norm{\psi}_{0,\dep}=\norm{\si^\delta\psi}_p,\ \
\norm{\psi}_{1,\dep}=\norm{\si^{1+\delta}\nabla\psi}_p+
\norm{\si^\delta \psi}_p,
\end{equation}
respectively, where $\norm{\psi}_p$ is  the $L^p$ norm.
\end{df}

\begin{lm}\label{thm:continuous}
For $p\geq 2$ and $0 < \delta < 2-3/p$, or $p=2$ and $\delta=-1$,
the operators $\na$ and $\vD$ are bounded linear maps from 
$W^{1,p}_\delta$ into $W^{0,p}_\delta$, and there is a continuous
embedding $W^{1,p}_{\delta}\subset W=W^{1,2}_{-1}$.
\end{lm}
\paragraph{Proof} It is obvious from the definitions that
the operators $\na$ and $\vD$ are bounded linear maps from 
$W^{1,p}_\delta$ into $W^{0,p}_{\delta+1}$. 
Suppose that
$p\geq 2$ and $0 <\delta < 2-3/p$. Then
$\si^{-(1+\delta)}$ is in $L_{\frac{2p}{p-2}}$. Set
$$
C=\, \norm{\si^{-(1+\delta)}}_{\frac{2p}{p-2}}
$$
which is a positive constant. We have
\begin{eqnarray*}
\norm{\psi}_{1,-1,2}&=&\norm{\na\psi}_2+ \norm{\si^{-1}\psi}_2\\
&=&\norm{\si^{-(1+\delta)}\si^{1+\delta}\na\psi}_2
+\norm{\si^{-(1+\delta)}\si^\delta\psi}_2\\
&\leq& \norm{\si^{-(1+\delta)}}_{\frac{2p}{p-2}}\norm{\na\psi}_p +
\norm{\si^{-(1+\delta)}}_{\frac{2p}{p-2}}\norm{\si^\delta \psi}_p\\
&=& C \norm{\psi}_{1,\delta,p}.\ \ \Box
\end{eqnarray*}

By Proposition \ref{thm:key},
$$
2\int_{M_r}|\vD\psi|^2 \geq
\frac{1}{10}\int_{M_r}|\na\psi|^2+\frac{1}{16}\int_\Om u^{-2}|du|^2|\psi|^2
+\sum_{\al=1}^\ell\int_{S^\al_r}\inn{\frac{H}{2}\psi-c(\nu)\vvD\psi}{\psi}
$$
for all $r\geq L$. So we have
\begin{equation}\label{eqn:Dd}
\int_M|\vD\psi|^2 \geq\frac{1}{20} \int_M|\na\psi|^2
\end{equation}
for $\psi\in C^\infty_0(N,S)$. By Lemma \ref{thm:continuous},
(\ref{eqn:Dd}) holds for $\psi\in W$, or equivalently
\begin{lm}\label{thm:Dbound}
For $\psi\in W$, we have
\begin{equation}\label{eqn:Dbound}
2\sqrt{5}\norm{\vD\psi}_2\, \geq\, \norm{\na\psi}_2.
\end{equation}
\end{lm} 

\begin{cor}
For $p\geq 2$ and $0 < \delta < 2-3/p$, or $p=2$ and $\delta=-1$,
the operator $\vD: W^{1,p}_\delta\to W^{0,p}_{\delta+1}$ is an injection.
\end{cor}
\paragraph{Proof} Suppose that  $\psi\in W^{1,p}_\delta$ and $\vD\psi=0$.
By Lemma \ref{thm:regularity}, $\psi$ is continuous on $M$ and smooth on
$M\setminus \pa\Om$. By Lemma \ref{thm:Dbound}, $\na\psi=0$, so 
$|\psi|$ is constant outside $M\setminus \pa\Om$. 
Then $|\psi|$ is constant on $M$ since $|\psi|$ is continuous on $M$.
We have
$$
\norm{\psi}_{1,\delta,p}\geq \norm{\si^\delta \psi}_p
=(\norm{ \si^{\delta p}  }_1)^{1/p} |\psi|
$$
where $\delta p\geq -1$, so $\norm{ \si^{\delta p} }_1=\infty$. We must
have $\psi\equiv 0$. $\ \ \Box$

\begin{lm}[ {\cite[Lemma 5.4 (a)]{Pa-Ta}}]\label{thm:zerod}
For sufficiently large $L$, we have
\begin{equation}\label{eqn:zerod}
\norm{\si^{-1}\psi}^2_{2;E^\al_{2L}}\leq 5 \norm{\na\psi}^2_{2;E^\al_{2L}}.
\end{equation}
for $\psi\in W$.
\end{lm}
 
\begin{lm}\label{thm:zeroD}
For sufficiently large $L$, there is a constant $c=c(L)>0$ such that for
all $\psi\in W$,
\begin{equation}\label{eqn:zeroD}
\norm{\psi}^2_W\leq c \norm{\vD \psi}^2_2.
\end{equation}
\end{lm}
\paragraph{Proof} 
We modify the proof of \cite[Lemma 5.5]{Pa-Ta}. 

Fix a large $L$ such that Lemma \ref{thm:zerod} holds. Choose a smooth 
function $\beta=\beta_L$ such that 
(i) $0\leq \beta\leq 1$, 
(ii) $\beta\equiv 1$ on each $E^\al_{3L}$,  
(iii) $\beta\equiv 0$ on $M_{2L}$, and
(iv) $|\na\beta|\leq 2/L$.
Given $\psi\in W$, let $\psi_\ii=(1-\beta)\psi$ and
$\psi_\ee=\beta\psi$. Note that the support of $\psi_\ii$ is contained
in $M_{3L}$, and the support of $\psi_\ee$ is
contained in $\cup_{\alpha=1}^\ell E^\alpha_{2L}$. 

We have
$$
\norm{\psi}^2_W=\norm{\psi}^2_{1,-1,2}
=(\norm{\na\psi}_2+\norm{\si^{-1}\psi}_2)^2\leq 2\norm{\na\psi}^2_2 +2 \norm{\si^{-1}\psi}^2_2.
$$
$$
\norm{\si^{-1}\psi}^2_2=\norm{\si^{-1}\psi}^2_{2;M_{2L}}+
\sum_{\al=1}^\ell \norm{\si^{-1}\psi}^2_{2;E^\al_{2L}}.
$$
Recall that $\si\geq 1$, so
$$
\norm{\si^{-1}\psi}^2_{2;M_{2L}}\leq \norm{\psi}^2_{2;M_{2L}}
=\norm{\psi_\ii}^2_{2;M_{2L}}\leq \norm{\psi_\ii}^2_2.
$$

If $\psi\in C^\infty(M,S)$ and  $\na\psi=0$, then $|\psi|$ is a
constant. So $B(\psi,\phi)= \int_M\inn{\na\psi}{\na\phi}$ is a positive
definite Hermitian form on $C^\infty_0(M,S)$. We have
\begin{equation}\label{eqn:interior}
\norm{\psi_\ii}^2_2\leq c_2(L)\norm{\na\psi_\ii}^2_2.
\end{equation}
(see e.g. \cite[Section 5.2]{Mo}). 
$$
\norm{\na\psi_\ii}_2^2=\norm{\na\psi-\na\psi_\ee}^2_2
\leq 2\norm{\na\psi}^2_2 +2 \norm{\na\psi_\ee}_2^2
$$
where
$$
\norm{\na\psi_\ee}^2_2=\norm{\na(\beta\psi)}^2_2\leq
2\norm{(\na\beta)\psi}^2_2+ 2\norm{\beta\na\psi}^2_2.
$$

Note that  $\na\beta$ vanishes outside $M_{3L}\setminus M_{2L}$, 
$|\na\beta|\leq 2/L$ and $1\leq 3L\si^{-1}$ on $M_{3L}\setminus M_{2L}$, so 
$$
\norm{(\na\beta)\psi}^2_2\leq
\norm{\frac{2}{L}\psi}^2_{2;M_{3L}\setminus M_{2L}}\leq 
\norm{\frac{2}{L}\cdot 3L\cdot\si^{-1}\psi}^2_{2;M_{3L}\setminus M_{2L}}
=36\norm{\si^{-1}\psi}^2_{2;M_{3L}\setminus M_{2L}}
$$
Recall that $|\beta|\leq 1$, so $\norm{\beta\na\psi}^2_2\leq \norm{\na\psi}^2_2$.
So we have
\begin{eqnarray*}
\norm{\si^{-1}\psi}^2_{2;M_{2L}}&\leq& \norm{\psi_\ii}^2_2
\leq c_2(L)\norm{\na\psi_\ii}^2_2
\leq 2c_2(L)(\norm{\na\psi}^2_2+\norm{\na\psi_\ee}^2_2)\\
&\leq&
2c_2(L)(3\norm{\na\psi}^2_2+72\norm{\si^{-1}\psi}^2_{2;M_{3L}\setminus
  M_{2L}})
\end{eqnarray*}
\begin{eqnarray*}
\norm{\si^{-1}\psi}^2_2&=&\norm{\si^{-1}\psi}^2_{2;M_{2L}}+
\sum_{\al=1}^\ell \norm{\si^{-1}\psi}^2_{2;E^\al_{2L}}\\
&\leq& 6c_2(L)\norm{\na\psi}^2_2 + (1+144 c_2(L))\left(\sum_{\al=1}^\ell
\norm{\si^{-1}\psi}^2_{2;E^\al_{2L}}\right)\\
&\leq& 6c_2(L)\norm{\na\psi}^2_2 + 5(1+144 c_2(L))\left(\sum_{\al=1}^\ell
\norm{\na\psi}^2_{2;E^\al_{2L}}\right)\\
&\leq&(726 c_2(L)+5) \norm{\na\psi}^2_2
\end{eqnarray*}
where we used Lemma \ref{thm:zerod}. So
$$
\norm{\psi}^2_W\leq 2\norm{\na\psi}^2_2+2\norm{\si^{-1}\psi}^2_2
\leq (1452c_2(L)+12)\norm{\na\psi}^2_2.
$$
By Lemma \ref{thm:Dbound}, $\norm{\na\psi}^2_2\leq 20
\norm{\vD\psi}^2_2$, so
$$
\norm{\psi}^2_W\leq 240(121c_2(L)+1)\norm{\vD\psi}^2_2. \ \ \Box
$$

\begin{lm}\label{thm:existW}
For each $\eta\in C^\infty_0(N,S)$ there exists a unique
$u\in W$ such that $\vD^2 u=\eta$. $\psi=\vD u\in W$
and
$$
\norm{\psi}^2_W\leq c\norm{\eta}^2_2.
$$
\end{lm}
\paragraph{Proof}
Consider the functional
$$
F(u)=\frac{1}{2}\norm{\vD u}^2_2 + \mathrm{Re}\inn{u}{\eta}_2
$$
on $W$. It is strictly convex, weakly lower semicontinuous, and
$$
F(u)\geq \frac{1}{2c}\norm{u}_W^2-\norm{u}_W\norm{\si\eta}_2,
$$
where $c$ is the constant in Lemma \ref{thm:zeroD}.
By calculus of variation, $F$ has a unique critical point $u\in W$
which is an absolute minimum in $W$, and $u$ is the weak solution to
$$
\vD^2 u=\eta.
$$
Let $\psi=\vD u\in W^{0,2}_0$. Then
$\psi$ is a weak solution to 
$$
\vD\psi=\eta.
$$
By Lemma \ref{thm:regularity}, $\psi\in W^{1,2}_\loc(M,S)\cap
C^\infty(M\setminus\pa\Om, S)$. By Lemma \ref{thm:zeroD},
$$
\norm{\psi}^2_W \leq c \norm{\eta}^2_2. 
$$
So $\psi\in W$. $\ \ \Box$

\begin{lm}\label{thm:aprioriI}
For each $\eta\in C^\infty_0(N,S)$ there exists a unique
$\psi\in W$ with $\vD\psi=\eta$. For $p, \delta$ as in Lemma
\ref{thm:continuous},
\begin{equation}\label{eqn:apriori}
\norm{\psi}_{1,\delta,p}\leq c(\delta,p)\norm{\eta}_{0,\delta+1,p}.
\end{equation}
\end{lm}
\paragraph{Proof} By Lemma \ref{thm:existW}, there exists $\psi\in W$
such that $\vD\psi=\eta$. The solution is unique by Lemma
\ref{thm:zeroD}.

By Lemma \ref{thm:regularity}, $\psi\in W^{1,p}_\loc(M,S)\cap
C^\infty(M\setminus\pa\Om, S)$. A slight modification of the proof of
\cite[Proposition 5.7]{Pa-Ta} gives the a priori bound
(\ref{eqn:apriori}).
$\ \ \Box$

\medskip

The proof of \cite[Proposition 5.9]{Pa-Ta} gives
\begin{pro}\label{thm:aprioriII}
Let $p,\delta$ be as in Lemma \ref{thm:continuous}. Then 
$\vD:W^{1,p}_\delta\to W^{0,p}_{\delta+1}$
is a surjection. If $\vD \psi=\eta$, where $\eta\in  W^{0,p}_\delta$
and $\psi\in W^{1,p}_\delta$, then the a priori bound (\ref{eqn:apriori}) 
holds.
\end{pro}

We now return to the proof of Theorem \ref{thm:dirac}.
Let $\psi_0$ be defined as in (\ref{eqn:ends}).
Then $\vD\psi_0\in L^2(M,S)\cap W^{0,p}_{1+\delta}\cap C^\infty(M,S)$
for all $2\leq p <\infty$ and $0<\delta<1-3/p$ (see \cite[Section
  4]{Pa-Ta} for details).
By Proposition \ref{thm:aprioriII}, there exists
$\psi_1\in W$ such that 
\begin{equation}
\vD\psi_1=-\vD\psi_0
\end{equation}
and $\psi_1\in W^{1,p}_\delta(M,S)$
for all $2\leq p<\infty$ and $0<\delta<1-3/p$. So
$$
\lim_{r\to\infty}r^{1-\ep}|\psi_1|=0
$$
for any $\ep>0$. By Lemma \ref{thm:regularity},
$\psi\in C^\infty(M,S)$ and $\psi\in W^{1,p}_\loc$ for any $2\leq
p<\infty$. Let $\psi=\psi_1+\psi_0$. Then $\psi$ is the desired
solution in Theorem \ref{thm:dirac}. The solution is unique by
Lemma \ref{thm:continuous}.

\subsection{Proof of Theorem \ref{thm:tgmass}}\label{sec:tgmass}

Let $\psi$ be the unique spinor given by Theorem \ref{thm:dirac}. 
By Proposition \ref{thm:key}, 
\begin{equation}
0\geq \frac{1}{10}\int_{M_r}|\na\psi|^2
+\frac{1}{16}\int_\Om u^{-2}|du|^2|\psi|^2
+\sum_{\al=1}^\ell\int_{S^\al_r}
\inn{\frac{H}{2}\psi-c(\nu)\vvD\psi}{\psi}.
\end{equation}
By calculations similar to those in \cite{Pa-Ta}, we have
$$
\lim_{r\to\infty}\int_{S^\al_r}\inn{\frac{H}{2}\psi-c(\nu)\vvD\psi}{\psi}
=-4\pi G m^\al_\infty |\psi^\al|^2
$$
where $m^\al_\infty$ is the ADM mass of the end $E^\al$. So
$$
4\pi G\sum_{\al=1}^\ell m^\al_\infty|\psi^\al|^2
\geq\frac{1}{10}\int_M|\na\psi|^2 +\frac{1}{16}\int_{\Om}u^{-2}|du|^2|\psi|^2.
$$

In particular, take $|\psi^\beta|=1$ and $\psi^\al=0$ for $\al\neq \beta$, we have
$$
4\pi G m^\beta_\infty
\geq \frac{1}{10}\int_M|\na\psi|^2
+\frac{1}{16}\int_\Om u^{-2}|du|^2|\psi|^2\geq 0.
$$
If $m^\beta_\infty=0$ for some $\beta$, we have $du=0$, so $u\equiv 1$, which
implies $\hg=\bg$. We also have $\na\psi=0$ on $M$.
We conclude that $\pa\Omega$ is connected, and $M$ is the Euclidean space.

\end{document}